\documentclass[12pt, twoside ]{article}
\usepackage{amsfonts}
\usepackage{amsmath,amsthm}
\usepackage{enumerate}
\usepackage{color}

\usepackage[OT1]{fontenc}
\usepackage{bm}
\usepackage{mathrsfs}
\allowdisplaybreaks[4]

\newtheorem{thm}{Theorem}[section]
\newtheorem{lem}[thm]{Lemma}

\numberwithin{equation}{section}

\begin{document}
\title{\textbf{Distribution of zeros and zero-density estimates for the derivatives of \textit{L}-functions attached to cusp forms}}
\author{Yoshikatsu Yashiro\\
\small Graduate School of Mathematics, Nagoya University,\\[-6pt] 
\small 464-8602 \ Chikusa-ku, Nagoya, Japan \\[-6pt] 
\small E-mail: m09050b@math.nagoya-u.ac.jp}
\date{}
\maketitle

\renewcommand{\thefootnote}{}

\vspace{-36pt}
\begin{abstract}
Let $f$ be a holomorphic cusp form of weight $k$ with respect to $SL_2(\mathbb{Z})$ which is a normalized Hecke eigenform, $L_f(s)$ the $L$-function attached to the form $f$. 
In this paper, we shall give the relation of the number of zeros of $L_f(s)$ and the derivatives of $L_f(s)$ using Berndt's method, and an estimate of zero-density of the derivatives of $L_f(s)$ based on Littlewood's method.
\end{abstract}

\section{Introduction}
Let $f$ be a cusp form of weight $k$  for $SL_2(\mathbb{Z})$ which is a normalized Hecke eigenform. Let $a_f(n)$ be the $n$-th Fourier coefficient of $f$ and set $\lambda_f(n)=a_f(n)/n^{(k-1)/2}$.  Rankin showed that  $\sum_{n\leq x}|\lambda_f(n)|^2=C_fx+O(x^{3/5})$ for $x\in\mathbb{R}_{>0}$, where $C_f$ is a positive constant depending on $f$ (see \cite[(4.2.3), p.364]{RAN}). 
The $L$-function attached to $f$ is defined by
\begin{align}
L_f(s)=\sum_{n=1}^\infty\frac{\lambda_f(n)}{n^s}=\prod_{p\text{:prime}}\left(1-\frac{\alpha_f(p)}{p^s}\right)^{-1}\left(1-\frac{\beta_f(p)}{p^s}\right)^{-1} \quad (\text{Re }s>1), \label{4LD}
\end{align}
where $\alpha_f(p)$ and $\beta_f(p)$ satisfy $\alpha_f(p)+\beta_f(p)=\lambda_f(p)$ and $\alpha_f(p)\beta_f(p)=1$. By Hecke's work (\cite{HEC}), the function $L_f(s)$ is analytically continued to the whole $s$-plane by
\begin{align}
(2\pi)^{-s-\frac{k-1}{2}}\Gamma(s+\tfrac{k-1}{2})L_f(s)=\int_0^\infty f(iy)y^{s+\frac{k-1}{2}-1}dy, \label{4AC}
\end{align}
and has a functional equation 
\begin{align*}
L_f(s)=\chi_f(s)L_f(1-s)
\end{align*}
where $\chi_f(s)$ is given by
\begin{align}
\chi_f(s)=&(-1)^{-\frac{k}{2}}(2\pi)^{2s-1}\frac{\Gamma(1-s+\frac{k-1}{2})}{\Gamma(s+\frac{k-1}{2})} \notag\\
=&2(2\pi)^{-2(1-s)}\Gamma(s+\tfrac{k-1}{2})\Gamma(s-\tfrac{k-1}{2})\cos\pi(1-s). \label{4XFE}
\end{align}
The second equality is deduced from the fact $\Gamma(s)\Gamma(1-s)=\pi/\sin(\pi s)$ and $\sin\pi(s+(k-1)/2)=(-1)^{k/2}\cos\pi(1-s)$. 
Similarly to the case of the Riemann zeta function $\zeta(s)$, it is conjectured that all complex zeros of $L_f(s)$ lie on the critical line  $\text{Re }s= 1/2$, namely, the Generalized Riemann Hypothesis (GRH). In order to support the truth of the GRH, the distribution and the density of complex zeros of $L_f(s)$ are studied without assuming the GRH. 

Lekkerkerker \cite{LEK} 
proved the approximate formula of a number of complex zeros of $L_f(s)$:
\begin{align}
N_f(T)=\frac{T}{\pi}\log\frac{T}{2\pi e}+O(\log T), \label{ML0}
\end{align}
where $T>0$ is sufficiently large, and $N_f(T)$ denotes the number of complex zeros of $L_f(s)$ in $0<{\rm Im\;}s\leq T$. The formula (\ref{ML0}) is an analogy of $N(T)$ which denotes the number of complex zeros of $\zeta(s)$ in $0<{\rm Im\;}s\leq T$. Riemann \cite{RIE}
showed that 
\begin{align}
N(T)=\frac{T}{2\pi}\log\frac{T}{2\pi e}+O(\log T). \label{Z0A}
\end{align}
(Later von Mangoldt \cite{VOM} proved (\ref{Z0A}) rigorously.)  In the Riemann zeta function, the zeros of derivative of $\zeta(s)$ have a connection with RH. Speiser \cite{SPE} showed that the Riemann Hypothesis (RH) is equivalent to the non-existence of complex zero of $\zeta'(s)$ in $\text{Re }s<1/2$, where $\zeta'(s)$ denotes the derivative function of $\zeta(s)$. Levinson and Montgomery \cite{L&M} proved that if RH is true, then $\zeta^{(m)}(s)$ has at most finitely many complex zeros in $0<\text{Re }s<1/2$ for any $m\in\mathbb{Z}_{\geq0}$. 

There are many studies of the zeros of $\zeta^{(m)}(s)$ without assuming RH. Spira \cite{SP1}, \cite{SP2} showed that there exist $\sigma_{m}\geq(7m+8)/4$ and $\alpha_{m}<0$ such that $\zeta^{(m)}(s)$ has no zero for ${\rm Re\;}s\leq\sigma_{m}$ and ${\rm Re\;}s\leq\alpha_m$, and exactly one real zero in each open interval $(-1-2n,1-2n)$ for $1-2n\leq\alpha_m$.
Later, Y{\i}ld{\i}r{\i}m \cite{YIL} showed that $\zeta''(s)$ and $\zeta'''(s)$ have no zeros in the strip $0\leq\text{Re\;}s<1/2$.  
Berndt \cite{BER} gave the relation of the number of complex zeros of $\zeta(s)$ and $\zeta^{(m)}(s)$:
\begin{align}
N_m(T)=N(T)-\frac{T\log 2}{2\pi}+O(\log T), \label{ZMA}
\end{align}
where $m\in\mathbb{Z}_{\geq1}$ is fixed and $N_m(T)$ denotes the number of complex zeros of $\zeta^{(m)}(s)$ in $0<{\rm Im\;}s\leq T$.  
Recently, Aoki and Minamide \cite{A&M} studied the density of zeros of $\zeta^{(m)}(s)$ in the right hand side of critical line ${\rm Re\;}s=1/2$ by using Littlewood's method. Let $N_m(\sigma,T)$ be the number of zeros of $\zeta^{(m)}(s)$ in $\text{Re\;}s\geq\sigma$ and $0<\text{Im\;}s\leq T$. They showed that 
\begin{align}
N_m(\sigma,T)=O\left(\frac{T}{\sigma-1/2}\log\frac{1}{\sigma-1/2}\right), \label{MME}
\end{align}
uniformly for $\sigma>1/2$. 
From (\ref{ZMA}) and (\ref{MME}), we see that almost all complex zeros of $\zeta^{(m)}(s)$ lie in the neighbourhood of the critical line. 

The purpose of this paper is to study the corresponding results of Berndt, Aoki and Minamide for the derivatives of $L_f(s)$, namely, the relation between the number of complex zeros of $L_f(s)$ and that of $L_f^{(m)}(s)$, and the density of zeros of $L_f^{(m)}{(s)}$ in the right half plane ${\rm Re\;}s>1/2$. 
Let $n_f$ be the smallest integer greater than 1 such that $\lambda_f(n_f)\ne0$.
Here $L^{(m)}_f(s)$ denotes the $m$-th derivative of $L_f(s)$ given by 
\begin{align*}
L_f^{(m)}(s)=\sum_{n=1}^\infty\frac{\lambda_f(n)(-\log n)^m}{n^s}=\sum_{n=n_f}^\infty\frac{\lambda_f(n)(-\log n)^m}{n^s} \quad (\text{Re }s>1),
\end{align*}
Differentiating both sides of (\ref{4AC}), we find that $L^{(m)}_f(s)$ is holomorphic in the whole $s$-plane and has the functional equation:
\begin{align}
L^{(m)}_f(s)=\sum_{r=0}^m\binom{m}{r}(-1)^{r}\chi_f^{(m-r)}(s)L_f^{(r)}(1-s).  \label{4DFE}
\end{align}

First in order to achieve the above purpose, we shall show the zero free regions for $L_f^{(m)}(s)$ by following Berndt's method (see \cite{BER}) and Spira's method (see \cite{SP1}, \cite{SP2}). 
\begin{thm}\label{THM0}
The following assertions hold for any $m\in\mathbb{Z}_{\geq0}$.
\begin{enumerate}[(i)]
\item\label{LRZF} There exists $\sigma_{f,m}\in\mathbb{R}_{>1}$ such that $L_f^{(m)}(s)$ has no zero for ${\rm Re\;}s\geq\sigma_{f,m}$.  
\item\label{LLZF} For any $\varepsilon\in\mathbb{R}_{>0}$, there exists $\delta_{f,m,\varepsilon}\in\mathbb{R}_{>(k-1)/2+1}$ such that $L_f^{(m)}(s)$ has no zero for $|s|\geq\delta_{f,m,\varepsilon}$ satisfying ${\rm Re\;}s\leq -\varepsilon$ and $|{\rm Im\;}s|\geq \varepsilon$.
\item\label{LOZF} There exists $\alpha_{f,m}\in\mathbb{R}_{<-(k-1)/2-1}$ such that $L_f^{(m)}(s)$ has only real zeros for ${\rm Re\;}s\leq\alpha_{f,m}$, and one real zero in each interval $(n-1,n)$ for $n\in\mathbb{Z}_{\leq \alpha_{f,m}}$.   
\end{enumerate}
\end{thm}

Next, based on Berndt's proof, we can obtain the following formula of the numbers of complex zeros of $L_f^{(m)}(s)$:
\begin{thm}\label{THM1}
For any fixed $m\in\mathbb{Z}_{\geq1}$, let $N_{f,m}(T)$ be the number of complex zeros of $L_f^{(m)}(s)$ in $0<{\rm Im\;}s\leq T$. 
Then for any large $T>0$, we have
\begin{align*}
N_{f,m}(T)=\frac{T}{\pi}\log\frac{T}{2\pi e}-\frac{T}{2\pi}\log{n_f}+O(\log T).
\end{align*}
Moreover the relation between $N_f(T)$ and $N_{f,m}(T)$ are given by
\begin{align*}
N_{f,m}(T)=N_{f}(T)-\frac{T}{2\pi}\log{n_f}+O(\log T).
\end{align*}

\end{thm}
Finally using the mean value formula for $L_f^{(m)}(s)$ obtained in \cite{YY3} and Littlewood's method, we obtain the estimate of density of zeros:
\begin{thm}\label{THM2}
For any $m\in\mathbb{Z}_{\geq0}$, let $N_{f,m}(\sigma,T)$ be the number of complex zeros of $L_f^{(m)}(s)$ in ${\rm Re\;}s\geq\sigma$ and $0<{\rm Im\;}s\leq T$. 
For any large $T>0$, we have
\begin{align}
N_{f,m}(\sigma,T)=O\left(\frac{T}{\sigma-1/2}\log\frac{1}{\sigma-1/2}\right) \label{ZD2}
\end{align}
uniformly for $1/2<\sigma\leq 1$. More precisely we have
\begin{align}
&\hspace{-12pt}N_{f,m}(\sigma,T)\notag\\
\leq&
\displaystyle\frac{2m+1}{2\pi}\frac{T}{\sigma-1/2}\log\frac{1}{\sigma-1/2}+\frac{1}{2\pi}\frac{T}{\sigma-1/2}\log\frac{(2m)!n_fC_f}{|\lambda_f(n_f)|^2(\log n_f)^{2m}}+ & \notag\\
&+O(\log T)+\frac{1}{2\pi}\frac{T}{\sigma-1/2}\times & \notag\\[9pt]
&\times\begin{cases} \displaystyle\log\left(1+O\left(\frac{(2\sigma-1)^{2m+1}(\log T)^{2m}}{T^{2\sigma-1}}\right)\right), & 1/2<\sigma<1, \\ \displaystyle\log\left(1+O\left(\frac{(2\sigma-1)^{2m+1}(\log T)^{2m+2}}{T}\right)\right), \hspace{6pt} & \sigma=1,\\ \displaystyle\log\left(1+O\left(\frac{(2\sigma-1)^{2m+1}}{T}\right)\right), & 1<\sigma<\sigma_{f,m}, \end{cases} 
\label{ZD1}
\end{align}
where $\sigma_{f,m}$ is given by (\ref{LRZF}) of Theorem \ref{THM0}.
\end{thm}

\subsection{Proof of Theorem \ref{THM0}}

First in order to show \eqref{LRZF}, we write 
$L_f^{(m)}(s)=\lambda_f(n_f)(-\log n_f)^m F(s)n_f^{-s}$ where 
\begin{align}
F(s)=1+\sum_{n=n_f+1}^\infty\frac{\lambda_f(n)}{\lambda_f(n_f)}\left(\frac{\log n}{\log n_f}\right)^m\left(\frac{n_f}{n}\right)^s \quad ({\rm Re\;}s>1). \label{FFF}
\end{align}
Deligne's result $|\lambda_f(n)|\leq d(n)\ll n^\varepsilon$ gives that there exist $c_f\in\mathbb{R}_{>0}$ and $\sigma_{f,m}\in\mathbb{R}_{>1}$ depending on $f$ and $m$ such that 
\begin{align}
|F(\sigma+it)-1|
\leq&\sum_{n=n_f+1}^\infty\left|\frac{\lambda_f(n)}{\lambda_f(n_f)}\right|\left(\frac{\log n}{\log n_f}\right)^m\left(\frac{n_f}{n}\right)^\sigma\notag\\
\leq& c_f\sum_{n=n_f+1}^\infty\frac{(\log{n}/\log{n_f})^m}{(n/n_f)^{\sigma-\varepsilon}}\leq\frac{1}{2}\label{FMZ}
\end{align}
for $\sigma\in\mathbb{R}_{\geq\sigma_{f,m}}$ and $t\in\mathbb{R}$, where $\varepsilon$ is an arbitrary positive number.
Hence $L_f^{(m)}(s)$ has no zeros for ${\rm Re\;}s\geq\sigma_{f,m}$, that is, \eqref{LRZF} is showed. 

Next we shall show \eqref{LLZF} and \eqref{LOZF}. 
Replacing $s$ to $1-s$ in \eqref{4XFE} and \eqref{4DFE}, we have  
\begin{align}
&(-1)^mL_f^{(m)}(1-s)\notag\\
&=\sum_{r=0}^m\binom{m}{r}L_f^{(m-r)}(s)\left(2(2\pi)^{-2s}\cos(\pi s)\Gamma(s-\tfrac{k-1}{2})\Gamma(s+\tfrac{k-1}{2})\right)^{(r)}. \label{CLZ1}
\end{align}
By the facts $(\cos(\pi s))^{(r)}=\pi^r(a_r\cos(\pi s)+b_r\sin(\pi s))$ where $a_r, b_r \in\{0,\pm 1\}$ and $((2\pi)^{-2s})^{(r)}=(-2\log 2\pi)^r(2\pi)^{-2s}$ for $r\in\mathbb{Z}_{\geq 0}$, the formula \eqref{CLZ1} is written as 
\begin{align}
(-1)^mL_f^{(m)}(1-s)
=2(2\pi)^{-2s}\sum_{r=0}^mR_{m-r}(s)(\Gamma(s-\tfrac{k-1}{2})\Gamma(s+\tfrac{k-1}{2}))^{(r)}, \label{CLZ2}
\end{align}
where
\begin{align*}
R_{m-r}(s)=\cos(\pi s)\sum_{j=0}^{m-r}a_j'L_f^{(j)}(s)+\sin(\pi s)\sum_{j=0}^{m-r}b_j'L_f^{(j)}(s)
\end{align*}
and $a_r', b_r'\in\mathbb{R}$. It is clear that $a_0'=1$, $b_0'=0$ and $R_0(s)=L_f(s)\cos(\pi s)$. Moreover we write \eqref{CLZ2} as
\begin{align}
\frac{(-1)^mL_f^{(m)}(1-s)}{2(2\pi)^{-2s}}=f(s)+g(s) \label{CLZ3}
\end{align} 
where
\begin{align*}
f(s)=&R_0(s)(\Gamma(s-\tfrac{k-1}{2})\Gamma(s+\tfrac{k-1}{2}))^{(m)},\\
g(s)=&\sum_{r=0}^{m-1}R_{m-r}(s)(\Gamma(s-\tfrac{k-1}{2})\Gamma(s+\tfrac{k-1}{2}))^{(r)}. 
\end{align*}
The formula \eqref{CLZ3} implies that if $|f(s)|>|g(s)|$ in a some region then $L_f^{(m)}(s)$ has no zero in this region. 

In order to investigate the behavior of $f(s)$ and $g(s)$, we shall consider the approximate formula for $(\Gamma(s-\tfrac{k-1}{2})\Gamma(s+\tfrac{k-1}{2}))^{(r)}/(\Gamma(s-\tfrac{k-1}{2})\Gamma(s+\tfrac{k-1}{2}))$. 
By Stirling's formula, it is known that 
\begin{align*}
\frac{\Gamma'}{\Gamma}(s)=\log s-\frac{1}{2s}+\int_0^\infty\frac{\{u\}-1/2}{(u+s)^2}du 
\end{align*}
for $s\in\mathbb{C}$ such that $|\arg s|\leq\pi-\delta$ where $\delta\in\mathbb{R}_{>0}$  is fixed (see \cite[Theorem A.5 b)]{K&V}). 
Writing the right-hand side of the above formula to $G^{(1)}(s)$ and putting $G^{(j)}(s)=(d^{j-1}/ds)G^{(1)}(s)$ for $j\in\mathbb{Z}_{\geq2}$, we shall use the following lemma:
\begin{lem}[{\cite[Lemma 2.3]{YY3}}] \label{LIB}
Let $F$ and $G$ be holomorphic function in the region $D$ such that $F(s)\ne0$ and $\log F(s)=G(s)$ for $s\in D$. Then for any fixed $r\in\mathbb{Z}_{\geq1}$, there exist $l_1,\cdots,l_r\in\mathbb{Z}_{\geq0}$ and $C_{(l_1,\cdots,l_r)}\in\mathbb{Z}_{\geq0}$ such that
\begin{align*}
\frac{F^{(r)}}{F}(s)
=\sum_{1l_1+\cdots+rl_r=r}C_{(l_1,\cdots,l_r)}(G^{(1)}(s))^{l_1}\cdots(G^{(r)}(s))^{l_r} 
\end{align*}
for $s\in D$. Especially $C_{(r,0,\cdots,0)}=1$. 
\end{lem}
The estimates 
\begin{align*}
|u+s|^2=&u^2+|s|^2+2u|s|\cos\arg s\\\geq&\begin{cases} |s|^2, & u\leq|s|, \; |\arg s|\leq\pi/2, \\ |s|^2(\sin\arg s)^2, & u\leq|s|,\; \pi/2\leq|\arg s|\leq\pi-\delta, \\ u^2, & u\geq|s|, \;|\arg s|\leq\pi/2, \\  4(1+\cos\arg s)u^2, & u\geq|s|,\; \pi/2\leq|\arg s|\leq\pi-\delta 
 \end{cases}
\end{align*}
give 
\begin{align*}
\int_{0}^{\infty}\frac{\{u\}-1/2}{(u+s)^{j+1}}du 
\ll&\int_0^{|s|}\frac{du}{|s|^{j+1}}+\int_{|s|}^\infty\frac{du}{u^{j+1}}\ll\frac{1}{|s|^{j}}
\end{align*}
for $|\arg s|\leq\pi-\delta$ and $j\in\mathbb{Z}_{\geq1}$. Then $G^{(j)}(s)$ is calculated as  
\begin{align*}
G^{(1)}(s)=&\log s+O\left(\frac{1}{|s|}\right),\\
G^{(j)}(s)=&\frac{(-1)^{j-1}(j-2)!}{s^{j-1}}+\frac{(-1)^j(j-1)!}{2s^j}+(-1)^{j+1}j!\int_0^\infty\frac{\{u\}-1/2}{(u+s)^{j+1}}du\\
=&O\left(\frac{1}{|s|^{j-1}}\right)
\end{align*}
for $j\in\mathbb{Z}_{\geq2}$. 
Hence the approximate formula for $(\Gamma^{(r)}/\Gamma)(s)$ is written as
\begin{align}
\frac{\Gamma^{(r)}}{\Gamma}(s)
=&(G^{(1)}(s))^r+O\left(\sum_{1q_1+\cdots+rq_r=r, \ q_1\ne r \;}\prod_{j=1}^r|G^{(j)}(s)|^{q_j}\right)\notag\\
=&\left(\log s+O\left(\frac{1}{|s|}\right)\right)^r+O\left(\frac{|\log s|^{r-1}}{|s|}\right)
=(\log s)^r\sum_{j=0}^r\frac{M_j(s)}{(\log s)^{j}} \label{CLZ4}
\end{align}
for $|\arg s|\leq\pi-\delta$ and $r\in\mathbb{Z}_{\geq0}$,  
where $M_j(s)$ is given by $M_j(s)=O(1/|s|^j)$ for $j\in\mathbb{Z}_{\geq 1}$ and $M_0(s)=1$. Using \eqref{CLZ4} and the trivial estimates $|1\pm({k-1})/{2s}|<2<|s|^{1/2}$,
\begin{align}
&|s\pm\tfrac{k-1}{2}|=|s|\times|1\pm\tfrac{k-1}{2s}|\gg|s|,\notag\\
&\log|s\pm\tfrac{k-1}{2}|=(\log|s|)\times\left(1+\frac{\log|1\pm\tfrac{k-1}{2s}|}{\log|s|}\right)\gg\log|s|, \label{CLE2}\\
&\log(s\pm\tfrac{k-1}{2})\ll\sqrt{(\log|s\pm\tfrac{k-1}{2}|)^2+(\arg(s\pm\tfrac{k-1}{2}))^2}\ll\log|s|\notag
\end{align}
for $|s|>(k-1)/2$, we obtain a desired formula:
\begin{align}
&\frac{(\Gamma(s-\tfrac{k-1}{2})\Gamma(s+\tfrac{k-1}{2}))^{(l)}}{\Gamma(s-\tfrac{k-1}{2})\Gamma(s+\tfrac{k-1}{2})}\notag\\
=&\sum_{j=0}^l\binom{l}{j}\frac{\Gamma^{(j)}}{\Gamma}(s-\tfrac{k-1}{2})\frac{\Gamma^{(l-j)}}{\Gamma}(s+\tfrac{k-1}{2})\notag\\
=&\sum_{j=0}^l\binom{l}{j}(\log(s-\tfrac{k-1}{2}))^j(\log(s+\tfrac{k-1}{2}))^{l-j}\times\notag\\
&\times\sum_{\begin{subarray}{c} 0\leq j_1+j_2\leq l,\\ 0\leq j_1\leq j,\; 0\leq j_2\leq l-j \end{subarray}}\frac{M_{j_1}(s-\tfrac{k-1}{2})}{(\log(s-\tfrac{k-1}{2}))^{j_1}}\frac{M_{j_2}(s+\tfrac{k-1}{2})}{(\log(s+\tfrac{k-1}{2}))^{j_2}}\notag\\
=&S_l(s)+T_l(s) \label{CLZ5}
\end{align}
for $l\in\mathbb{Z}_{\geq0}$ and $s\in\mathbb{C}$ such that $|s|>(k-1)/2$ and $|\arg s|\leq\pi-\delta$ , where $S_l(s)$, $T_l(s)$ are given by
\begin{align*}
S_l(s)=&(\log(s-\tfrac{k-1}{2})+\log(s+\tfrac{k-1}{2}))^l,\\
T_l(s)=&O\left(\frac{1}{|s|\log|s|}\sum_{j=0}^l(\log|s|)^j(\log|s|)^{l-j}\right)=O\left(\frac{(\log |s|)^{l-1}}{|s|}\right)
\end{align*}
respectively for $l\in\mathbb{Z}_{\geq 1}$, especially $S_0(s)=1$ and $T_0(s)=0$.

Next using $R_r(s)$, $S_r(s)$ and $T_r(s)$, we shall write a condition which gives $|f(s)|>|g(s)|$ for some region. From \eqref{CLZ3} and \eqref{CLZ5}, the inequality $|f(s)|>|g(s)|$ is equivalent to
\begin{align*}
|S_m(s)+T_m(s)|>\left|\sum_{r=0}^{m-1}\frac{R_r}{R_0}(s)(S_r(s)+T_r(s))\right|. 
\end{align*}
Dividing the both sides of the above formula by $S_{m-1}(s)$ and applying the triangle inequality, we see that if 
\begin{align}
|S_1(s)|>\left|\frac{T_m}{S_{m-1}}(s)\right|+\left|\sum_{r=0}^{m-1}\frac{R_r}{R_0}(s)\left(\frac{1}{S_{m-1-r}}(s)+\frac{T_r}{S_{m-1}}(s)\right)\right| \label{CLZ6}
\end{align}
is true, then $|f(s)|>|g(s)|$ is true for $|s|>(k-1)/2$ and $|\arg s|\leq\pi-\delta$. To show the truth of \eqref{CLZ6}, we shall consider an upper bound of $(1/S_r)(s)$, $T_r(s)$ and $(R_r/R_0)(s)$. 
The estimates \eqref{CLE2} and $|\log z|\geq\log|z|$ for $z\in\mathbb{C}$ give
\begin{align}
\left|\frac{1}{S_{r}}(s)\right|\leq\frac{1}{(\log|s-\tfrac{k-1}{2}|+\log|s+\tfrac{k-1}{2}|)^r}\leq\frac{C_1}{(\log|s|)^r} \label{CLZ7}
\end{align}
for the above $s$ and $r\in\mathbb{Z}_{\geq 0}$, here and later $C_1, C_2, \dots$ denote positive constants depending on $f$, $r$ and $\delta$.  
Since $L_f^{(j)}(s)$ and $(1/L_f)(s)$  are absolutely convergent for ${\rm Re\;}s>1$, it follows that 
\begin{align}
\left|\frac{R_r}{R_0}(s)\right|=\left|\sum_{j=0}^ra'_j\frac{L_f^{(j)}}{L_f}(s)+\tan(\pi s)\sum_{j=0}^rb'_j\frac{L_f^{(j)}}{L_f}(s)\right|\leq C_2+C_3|\tan(\pi s)| \label{CLZ9}
\end{align}
for ${\rm Re\;}s\geq1+\varepsilon$. Here $\tan(\pi s)$ is estimated as
\begin{align}
|\tan\pi(\sigma+it)|=\left|\frac{e^{-2t}e^{2\pi i\sigma}-1}{e^{-2t}e^{2\pi i\sigma}+1}\right|\leq\begin{cases} 2/(1-e^{-2\varepsilon}), & \text{if } |t|\geq\varepsilon, \\ 3, & \text{if } \sigma\in\mathbb{Z} \end{cases} \label{CLZ10}
\end{align}
where $\varepsilon$ is a fixed positive number.  
   
Combining \eqref{CLZ7}--\eqref{CLZ10}, we see that the right-hand side of \eqref{CLZ6} is estimated as 
\begin{align}
&\left|\frac{T_m}{S_{m-1}}(s)\right|+\left|\sum_{r=0}^{m-1}\frac{R_r}{R_0}(s)\left(\frac{1}{S_{m-1-r}}(s)+\frac{T_r}{S_{m-1}}(s)\right)\right|\notag\\
&\leq\frac{C_4}{|s|}+C_5|\tan(\pi s)|\sum_{r=0}^{m-1}\left(\frac{1}{(\log|s|)^{m-1-r}}+\frac{1}{|s|(\log|s|)^{m-r}}\right)\leq C_{f,m,\delta,\varepsilon} \label{CLZ11} 
\end{align}
for $|s|>(k-1)/2$ and ${\rm Re\;}s\geq1+\varepsilon$ provided $|{\rm Im\;}s|\geq\varepsilon$ or ${\rm Re\;}s\in\mathbb{Z}$, where $C_{f,m,\delta,\varepsilon}$ is a positive constant depending on $f$, $m$, $\delta$ and $\varepsilon$. Fix $\delta=\tan^{-1}(2\varepsilon/(k-1))$ and choose $r_{f,m}\in\mathbb{R}_{>(k-1)/2}$ such that $C_{f,m,\delta,\varepsilon}<(\log r_{f,m})/C_1$. The inequalities \eqref{CLZ7} and \eqref{CLZ11} imply that \eqref{CLZ6} is true, that is, $L_f^{(m)}(1-s)$ has no zero for $s\in\mathbb{C}$ such that $|s|\geq r_{f,m}$, ${\rm Re\;}s\geq1+\varepsilon$ and $|{\rm Im\;}s|\geq\varepsilon$. Therefore, we conclude that for any $\varepsilon\in\mathbb{R}_{>0}$ there exists $\delta_{f,m,\varepsilon}\in\mathbb{R}_{>(k-1)/2+1}$ such that $L_f^{(m)}(s)$ has no zero in the region $|s|\geq\delta_{f,m,\varepsilon}$, ${\rm Re\;}s\leq-\varepsilon$ and $|{\rm Im\;}s|>\varepsilon$, that is, the proof of \eqref{LLZF} is completed.

Finally we shall show \eqref{LOZF} applying Rouch\'{e}'s theorem to $f(s)$ and $g(s)$. For $n\in\mathbb{Z}_{\geq1}$ let $D_n$ be the region $n\leq{\rm Re\;}s\leq n+1$ and $|{\rm Im\;}s|\leq 1/2$. By \eqref{CLZ10} and \eqref{CLZ11}, we see that there exists $\delta_{f,m,1/2}\in\mathbb{R}_{>(k-1)/2}$ such that $|f(s)|>|g(s)|$ is true in the boundary of $D_n$ and the region $|s|>\delta_{f,m,1/2}$ and ${\rm Re\;}s\geq1+1/2$. Then the number of zeros of $f(s)$ is equal to that of $f(s)+g(s)$ in the interior of $D_n$. From \eqref{CLZ2}, \eqref{CLZ3} and \eqref{CLZ5}, the function $f(s)$ is written as
\begin{align}
f(s)=L_f(s)\cos(\pi s)\Gamma(s-\tfrac{k-1}{2})\Gamma(s+\tfrac{k-1}{2})(S_m(s)+T_m(s)). \label{FST}
\end{align}
When $R_{f,m}$ is chosen such that $R_{f,m}\geq\delta_{f,m,1/2}$ and $C_6/(R_{f,m}\log R_{f,m})<1$, 
the formula \eqref{CLZ7} gives
\begin{align}
\left|\frac{T_m}{S_m}(s)\right|\leq \frac{C_6}{|s|\log|s|}<1, \label{TS1}
\end{align}
that is, $S_m(s)+T_m(s)$ has no zero for $|s|\geq R_{f,m}$. Hence $f(s)$ has the only real zero $s=n+1/2$ in $D_n$. It is clear to show 
$\overline{f(s)}=f(\overline{s})$ and $\overline{g(s)}=g(\overline{s})$ for $s\in\mathbb{C}$, which imply that $L_f^{(m)}(1-s)$ has a only real zero in the interior of $D_n$. Replacing $1-s$ to $s$, we obtain the fact that there exists $\alpha_{f,m}\in\mathbb{R}_{<-(k-1)/2-1}$ such that $L_f^{(m)}(s)$ has no complex zero for ${\rm Re\;}s<\alpha_{f,m}$ and one real zero in each open interval $(n-1,n)$ for $n\in\mathbb{Z}_{\leq\alpha_{f,m}}$. The proof of \eqref{LOZF} is completed.   

\subsection{Proof of Theorem \ref{THM1}}
Using Theorem \ref{THM0}, we can choose $\alpha_{f,m}\in\mathbb{R}_{<-(k-1)/2}$ and $\sigma_{f,m}\in\mathbb{R}_{>1}$ such that $L_f^{(m)}(s)$ has no zeros in the region ${\rm Re\;}s\leq\alpha_{f,m}$ and ${\rm Re\;}s\geq\sigma_{f,m}$. 
Moreover, choose $\tau_{f,m}\in\mathbb{R}_{>2}$ and $T\in\mathbb{R}_{>0}$ such that $L_f^{(m)}(s)$ has no zeros for $0<{\rm Im\;}s\leq\tau_{f,m}$ and ${\rm Im\;}s=T$. Using the residue theorem in the region $\alpha_{f,m}\leq{\rm Re\;}s\leq\sigma_{f,m}$ and $\tau_{f,m}\leq{\rm Im\;}s\leq T$, we get
\begin{align}
N_{f,m}(T)
=&\frac{1}{2\pi i}\left(\int_{\alpha_{f,m}+i\tau_{f,m}}^{\sigma_{f,m}+i\tau_{f,m}}+\int_{\sigma_{f,m}+i\tau_{f,m}}^{\sigma_{f,m}+iT}+\int_{\sigma_{f,m}+iT}^{\alpha_{f,m}+iT}+\right.\notag\\
&\left.+\int_{\alpha_{f,m}+iT}^{\alpha_{f,m}+i\tau_{f,m}}\right)(\log L_f^{(m)}(s))'ds
=:I_1+I_2+I_3+I_4. \label{LDB}
\end{align}
First, it is clear that
\begin{align}
I_1=\frac{\log L_f^{(m)}(\sigma_{f,m}+i\tau_{f,m})-\log L_f^{(m)}(\alpha_{f,m}+i\tau_{f,m})}{2\pi i}
=O(1). \label{LDC}
\end{align}
To approximate $I_2$, we write $L_f^{(m)}(s)=\lambda_f(n_f)(-\log{n_f})^mF(s)n_f^{-s}$ where $F(s)$ is given by \eqref{FFF}. 
Using (\ref{FMZ}) we find that $1/2\leq|F(s)|\leq3/2$, ${\rm Re\;}F(s)\geq1/2$ and $|\arg F(s)|<\pi/2$ for $s=\sigma_{f,m}+it\;(t\in\mathbb{R})$. Hence $I_2$ is approximated as 
\begin{align}
I_2=&\frac{1}{2\pi i}\left[\log\frac{\lambda_f(n_f)(-\log{n_f})^m}{n_f^s}+\log F(s)\right]_{\sigma_{f,m}+i\tau_{f,m}}^{\sigma_{f,m}+iT}\notag\\
=&\frac{-(\sigma_{f,m}+iT)\log n_f}{2\pi i}+O(1)
=-\frac{T}{2\pi}\log{n_f}+O(1). \label{LDE}
\end{align}

Next we shall estimate $I_3$. The formula (\ref{4DFE}), the approximate functional equation for $L_f^{(m)}(s)$ (see \cite[Theorem 1.2]{YY3}) and Rankin's result, there exists $A\in\mathbb{R}_{\geq0}$ such that $L_f^{(m)}(\sigma+it)=O(|t|^A)$ uniformly for $\sigma\in[\alpha_{f,m},\sigma_{f,m}]$. It implies that
\begin{align}
I_3=&\frac{\log L_f^{(m)}(\alpha_{f,m}+iT)-\log L_f^{(m)}(\sigma_{f,m}+iT)}{2\pi i}\notag\\
=&\frac{\arg L_f^{(m)}(\alpha_{f,m}+iT)-\arg L_f^{(m)}(\sigma_{f,m}+iT)}{2\pi}+O(\log T). \label{LDF}
\end{align}
To estimate the first term of the right-hand side of (\ref{LDF}), we write 
$L_f^{(m)}(\sigma+iT)=(-1)^me^{-iT\log n_f}\lambda_f(n_f)G(\sigma+iT)$ where 
\begin{align*}
G(\sigma+iT)=\frac{(\log{n_f})^m}{n_f^{\sigma}}+\frac{1}{\lambda_f(n_f)}\sum_{n=n_f+1}^\infty\frac{\lambda_f(n)(\log n)^m}{n^{\sigma}}e^{iT\log\frac{n_f}{n}}  
\end{align*}
for $\sigma\in\mathbb{R}_{>1}$. Let $Q$ be the number of zeros of ${\rm Re\;}G(s)$ on the line segment $(\alpha_{f,m}+iT,\sigma_{f,m}+iT)$. Divide this line into $Q+1$ subintervals by these zeros.
 Then the sign of ${\rm Re\;}G(s)$ is constant, and the variation of $\arg G(s)$ is at most $\pi$ on each subinterval. 
Hence, there exists constant $C$ such that $\arg G(s)=\arg L_f^{(m)}(s)+C$ on the divided line, it follows that
\begin{align}
|\arg L_f^{(m)}(\alpha_{f,m}+iT)-\arg L_f^{(m)}(\sigma_{f,m}+iT)|\leq (Q+1)\pi. \label{LDH}
\end{align}
In order to estimate $Q$, let $H(z)=(G(z+iT)+\overline{G(\overline{z}+iT)})/2$. Then we find that  
\begin{align}
H(\sigma)=&{\rm Re\;}G(\sigma+iT)\notag\\
=&\frac{(\log n_f)^m}{n_f^\sigma}\left(1+\sum_{n=n_f+1}^\infty\frac{\lambda_f(n)}{\lambda_f(n_f)}\left(\frac{\log n}{\log n_f}\right)^m\left(\frac{n_f}{n}\right)^\sigma\cos\left(T\log\frac{n_f}{n}\right)\right) \label{LDG}
\end{align}
for $\sigma\in\mathbb{R}_{>1}$. The formulas (\ref{FMZ}) and (\ref{LDG}) give 
\begin{align}
\frac{1}{2}\frac{(\log n_f)^m}{{n_f}^{\sigma_{f,m}}}\leq H(\sigma_{f,m})\leq\frac{3}{2}\frac{(\log n_f)^m}{{n_f}^{\sigma_{f,m}}}. \label{HEE}
\end{align}
Take $T$ sufficiently large such that $T-\tau_{f,m}>2(\sigma_{f,m}-\alpha_{f,m})$ if necessary. Since ${\rm Im}(z+iT)\geq T-(T-\tau_{f,m})>0$ for $z\in\mathbb{C}$ such that $|z-\sigma_{f,m}|<T-\tau_{f,m}$, it follows that $H(z)$ is analytic in the circle $|z-\sigma_{f,m}|<T-\tau_{f,m}$. Note that there exists a positive constant $B$ such that $H(z)=O(T^B)$ in this circle because of the fact that $L_f(\sigma+it)=O(|t|^A)$.
For $u\in\mathbb{R}_{\geq0}$, let $P(u)$ be the number of zeros of $H(z)$ in $|z-\sigma_{f,m}|\leq u$.
Then using the trivial estimate 
\begin{align*}
P(\sigma_{f,m}-\alpha_{f,m})&\leq\frac{1}{\log 2}\int_{\sigma_{f,m}-\alpha_{f,m}}^{2(\sigma_{f,m}-\alpha_{f,m})}\frac{P(u)}{u}du\\
&\leq\frac{1}{\log 2}\int_{0}^{2(\sigma_{f,m}-\alpha_{f,m})}\frac{P(u)}{u}du,
\end{align*}
Jensen's formula (see \cite[Chapter 3.61]{TID}), the above note and (\ref{HEE}), we have
\begin{align*}
&P(\sigma_{f,m}-\alpha_{f,m})\\
&\ll\int_0^{2(\sigma_{f,m}-\alpha_{f,m})}\frac{P(u)}{u}du\\
&=\frac{1}{2\pi}\int_0^{2\pi}\log|H(\sigma_{f,m}+2(\sigma_{f,m}-\alpha_{f,m})e^{i\theta})|d\theta-\log|H(\sigma_{f,m})|\\
&\ll\int_0^{2\pi}\log T^Bd\theta+1\ll \log T,
\end{align*}
Therefore $Q$ is estimated as
\begin{align}
Q=\#\{\sigma\in(\alpha_{f,m},\sigma_{f,m})\mid F(\sigma)=0\}\ll P(\sigma_{f,m}-\alpha_{f,m})\ll \log T. \label{LDI}
\end{align}
Combining (\ref{LDF}), (\ref{LDH}), (\ref{LDI}), we obtain the estimate of $I_3$:
\begin{align}
I_3=O(\log T). \label{LD3}
\end{align}

Finally in order to approximate $I_4$, we shall obtain the approximate formula for $\log L_f^{(m)}(\alpha_{f,m}+iT)$ as $T\to\infty$. 
By the proof of Theorem \ref{THM0}, there exists $\delta_{f,m}\in\mathbb{R}_{>0}$ such that 
\begin{align}
\left|\frac{g}{f}(1-s)\right|< 1, \quad \left|\frac{T_m}{S_m}(1-s)\right|<1 \label{FGST}
\end{align}
for $s\in\mathbb{C}$ in the region $|s-(1-(k-1)/2)|>\delta_{f,m}$, ${\rm Re\;}s<1-(k-1)/2$ and $|{\rm Im\;}s|>1/2$. Here choose $\alpha_{f,m}\in\mathbb{R}_{<0}$ such that $\alpha_{f,m}<1-(k-1)/2-\delta_{f,m}$ if necessary. Then the path of $I_4$ is contained in the above region. Replacing $s$ to $1-s$ and taking logarithmic function in the both sides of \eqref{CLZ3}, we obtain
\begin{align}
\log L_f^{(m)}(\alpha_{f,m}+iT)
=&-2(1-\alpha_{f,m}-iT)\log 2\pi+\log f(1-\alpha_{f,m}-iT)+\notag\\
&+\log\left(1+\frac{g}{f}(1-\alpha_{f,m}-iT)\right)+O(1). \label{I4A}
\end{align}
The first formula of \eqref{FGST} gives 
 $|\arg(1+(g/f)(1-\alpha_{f,m}-iT))|<\pi/2$ and 
\begin{align}
&\log\left(1+\frac{g}{f}(1-\alpha_{f,m}-iT)\right)\notag\\
&\ll\sqrt{\left|1+\frac{g}{f}(1-\alpha_{f,m}-iT)\right|^2+\left(\arg\left(1+\frac{g}{f}(1-\alpha_{f,m}-iT)\right)\right)^2}\ll 1. \label{LAFG}
\end{align}
By \eqref{FST}, the second term of the right-hand side of \eqref{I4A} is written as
\begin{align}
&\hspace{-12pt}\log f(1-\alpha_{f,m}-iT)\notag\\
=&\log\Gamma(1-\alpha_{f,m}-\tfrac{k-1}{2}-iT)+\log\Gamma(1-\alpha_{f,m}+\tfrac{k-1}{2}-iT)+\notag\\
&+\log S_m(1-\alpha_{f,m}-iT)+\log\left(1+\frac{T_m}{S_m}(1-\alpha_{f,m}-iT)\right)+\notag\\
&+\log L_f(1-\alpha_{f,m}-iT)+\log\cos\pi(1-\alpha_{f,m}-iT). \label{I4B}
\end{align}
Now it is clear that
\begin{align*}
&\cos\pi(1-\alpha_{f,m}-iT)=e^{\pi T}e^{i(1-\alpha_{f,m})-\log 2}(1+{e^{-2\pi(1-\alpha_{f,m})i}}/{e^{2\pi T}}),\\
&\log L_f(1-\alpha_{f,m}-iT)=\sum_{n=1}^\infty\frac{b_f(n)}{n^{1-\alpha_{f,m}-iT}}
\end{align*}
where $b_f(n)$ is given by
\begin{align*}
b_f(n)=\begin{cases} (\alpha_f(p)^r+\beta_f(p)^r)/r, & n=p^r, \\ 0 , & \text{otherwise} \end{cases}
\end{align*} 
and $\alpha_f(p)$, $\beta_f(p)$ are given by \eqref{4LD}. Hence the fifth and sixth terms of the right-hand sides of \eqref{I4B} are approximated as
\begin{align}
\log L_f(1-\alpha_{f,m}-iT)+\log\cos\pi(1-\alpha_{f,m}-iT)=\pi T+O(1).
\end{align}
By the similar discussion of \eqref{LAFG}, the fourth term of the right-hand sides of \eqref{I4B} is estimated as
\begin{align}
\log\left(1+\frac{T_m}{S_m}(1-\alpha_{f,m}-iT)\right)\ll 1.
\end{align}
The trivial approximate formula
\begin{align*}
\log(1-\alpha_{f,m}\pm\tfrac{k-1}{2}-iT)=\log T-(\pi/2)i+O(1/T)
\end{align*} 
gives that the third term of the right-hand sides of \eqref{I4B} is approximated as
\begin{align}
&\log S_m(1-\alpha_{f,m}-iT)\notag\\
&=m\log\left(\log(1-\alpha_{f,m}-\tfrac{k-1}{2}-iT)+\log(1-\alpha_{f,m}+\tfrac{k-1}{2}-iT)\right)\notag\\
&=m\log\log T+O(1).
\end{align}
Using Stirling's formula $$\log\Gamma(s)=(s-1/2)\log s-s+\log\sqrt{2\pi}+O(1/|s|)$$ and the approximate formula of $\log(1-\alpha_{f,m}\pm({k-1})/{2}-iT)$, we approximate the first and second terms of the right-hand sides of \eqref{I4B} as
\begin{align}
&\log\Gamma(1-\alpha_{f,m}+\tfrac{k-1}{2}-iT)+\log\Gamma(1-\alpha_{f,m}-\tfrac{k-1}{2}-iT)\notag\\
&=(1-2\alpha_{f,m}-2iT)(\log T-(\pi/2)i+O(1/T))-2(1-\alpha_{f,m}-iT)+O(1)\notag\\
&=-2iT\log(T/e)-\pi T+(1-2\alpha_{f,m})\log T+O(1). \label{I4B12}
\end{align}  
Combining (\ref{I4A})--(\ref{I4B12}), we obtain a desired approximate formula as $$\log L_f^{(m)}(\alpha_{f,m}+iT)=-2iT\log\frac{T}{2\pi e}+O(\log T),$$ which implies that
\begin{align}
I_4=\frac{T}{\pi}\log\frac{T}{2\pi e}+O(\log T). \label{LD4}
\end{align}
From (\ref{LDC}), (\ref{LDE}), (\ref{LD3}) and (\ref{LD4}), the proof of Theorem \ref{THM1} is completed.


\subsection{Proof of Theorem \ref{THM2}}
Write $L_f^{(m)}(s)=\lambda_f(n_f)(-\log n_f)^mF(s)/n_f^s$ where $F(s)$ is given by (\ref{FFF}). By the proof of Theorem \ref{THM0}, we can choose $\sigma_{f,m}\in\mathbb{R}_{>1}$ such that $L_f(s)$ has no zero for ${\rm Re\;}s>\sigma_{f,m}$ and
\begin{align*}
\sum_{n=n_f+1}^\infty\left|\frac{\lambda_f(n)}{\lambda_f(n_f)}\right|\left(\frac{\log n}{\log n_f}\right)^m\left(\frac{n_f}{n}\right)^{\sigma_{f,m}/2}\leq\frac{1}{2}.
\end{align*}
Note that \eqref{FMZ} and the above inequality give
\begin{align}
|F(s)-1|&
\leq\sum_{n=n_f+1}^\infty\left|\frac{\lambda_f(n)}{\lambda_f(n_f)}\right|\left(\frac{\log n}{\log n_f}\right)^m\left(\frac{n_f}{n}\right)^{\sigma_{f,m}/2+\sigma/2}
\leq\frac{1}{2}\left(\frac{n_f}{n_f+1}\right)^{\sigma/2}. \label{LAS}
\end{align}
for ${\rm Re\;}s\geq\sigma_{f,m}$. 
Applying Littlewood's formula (see \cite[chapter 3.8]{TID}) to $F(s)$, we obtain
\begin{align}
2\pi\sum_{\begin{subarray}{c}F(\rho)=0, \\ \sigma\leq{\rm Re\;}\rho\leq\sigma_{f,m},\\ 1\leq{\rm Im\;}\rho\leq T \end{subarray}}({\rm Re\;}\rho-\sigma)
=&\int_1^T\log|F(\sigma+it)|dt-\int_1^T\log|F(\sigma_{f,m}+it)|dt+\notag\\[-24pt]
&+\int_\sigma^{\sigma_{f,m}}\arg{F(u+iT)}dt-\int_\sigma^{\sigma_{f,m}}\arg{F(u+i)}dt\notag\\
=&:I_1+I_2+I_3+I_4 \label{LA1}
\end{align}
for $\sigma\in\mathbb{R}_{>1/2}$. Here we shall estimate $I_2$. Cauchy's theorem gives that
\begin{align}
I_2=\int_1^T\log|F(v+it)|dt+\int_{\sigma_{f,m}}^{v}\log|F(u+i)|du-\int_{\sigma_{f,m}}^{v}\log |F(u+iT)|du\label{L2A}
\end{align}
for all $v>\sigma_{f,m}$. The fact that $\log{|X|}\leq |X-1|$ for $X\in\mathbb{C}$ and $-\log|Y|\leq 2|Y-1|$ for $Y\in\mathbb{C}$ satisfying $|Y|\geq1/2$, and (\ref{LAS}) imply that
\begin{align}
\int_1^T\log|F(v+it)|dt\leq\frac{(T-1)}{2}\left(\frac{n_f}{n_f+1}\right)^{{v}/{2}}
\label{L2B}
\end{align}
and 
\begin{align}
&\int_{\sigma_{f,m}}^{v}\log|F(u+i)|du-\int_{\sigma_{f,m}}^{v}\log |F(u+iT)|du\notag\\
&\ll\int_{\sigma_{f,m}}^{v}\left(\frac{n_f}{n_f+1}\right)^{{u}/{2}}du\ll 1. \label{L2C}
\end{align}
Combining (\ref{L2A})--(\ref{L2C})  we get
\begin{align}
I_2=O(1). \label{LI2}
\end{align}
By the same discussion of an estimate of $I_3$ in proof of Theorem \ref{THM1}, we can obtain
\begin{align}
I_3+I_4=O(\log T). \label{LA3}
\end{align}
To estimate $I_1$, we calculate
\begin{align}
I_1=\frac{T-1}{2}\log\frac{n_f^{2\sigma}}{|\lambda_f(n_f)|^2(\log n_f)^{2m}}+\frac{1}{2}\int_1^T\log|L_f^{(m)}(\sigma+it)|^2dt. \label{LA4}
\end{align}
Jensen's inequality gives
\begin{align}
\int_1^T\log|L_f^{(m)}(\sigma+it)|^2dt\leq(T-1)\log\left(\frac{1}{T-1}\int_1^T|L_f^{(m)}(\sigma+it)|^2dt\right). \label{LA5}
\end{align}
Combining (\ref{LA1}), (\ref{LI2})--(\ref{LA5}), we obtain
\begin{align}
\sum_{\begin{subarray}{c}F(\rho)=0, \\ \sigma\leq{\rm Re\;}\rho\leq\sigma_{f,m},\\ 1\leq{\rm Im\;}\rho\leq T \end{subarray}}({\rm Re\;}\rho-\sigma)\leq&\frac{T-1}{4\pi}\log\left(\frac{1}{T-1}\int_1^T|L_f^{(m)}(\sigma+it)|^2dt\right)+\notag\\[-24pt]
&+\frac{T-1}{4\pi}\log\frac{n_f^{2\sigma}}{|\lambda_f(n_f)|^2(\log n_f)^{2m}}+O(\log{T}). \label{LAA}
\end{align}

First, we consider the mean square of  $L_f^{(m)}(s)$ for ${\rm Re\;}s>1$. Then we can calculate as follows:
\begin{align}
&\hspace{-12pt}\int_1^T|L_f^{(m)}(\sigma+it)|^2dt\notag\\
=&\sum_{n_1, n_2=1}^\infty\frac{\overline{\lambda_f(n_1)}\lambda_f(n_2)(\log n_1)^m(\log n_2)^m}{(n_1n_2)^\sigma}\int_{\max\{n_1,n_2\}}^T
\left(\frac{n_1}{n_2}\right)^{it}dt\notag\\
=&(T-1)\sum_{n=1}^\infty\frac{|\lambda_f(n)|^2(\log n)^{2m}}{n^{2\sigma}}+\notag\\
&+\frac{1}{i}\sum_{\begin{subarray}{c} n_1, n_2=1,\\ n_1\ne n_2 \end{subarray}}^\infty\frac{\overline{\lambda_f(n_1)n_1^{-iT}}\lambda_f(n_2)n_2^{-iT}(\log n_1)^m(\log n_2)^m}{(n_1n_2)^\sigma\log(n_1/n_2)}-\notag\\
&-\frac{1}{i}\sum_{\begin{subarray}{c} n_1, n_2=1,\\ n_1\ne n_2 \end{subarray}}^\infty\frac{\overline{\lambda_f(n_1)n_1^{-i\max\{n_1,n_2\}}}\lambda_f(n_2)n_2^{-i\max\{n_1,n_2\}}(\log n_1)^m(\log n_2)^m}{(n_1n_2)^\sigma\log(n_1/n_2)}. \label{L2I}
\end{align}
By the same discussion for $U_\sigma(x)$ with
\begin{align*}
(\alpha_{n_1},\beta_{n_2})=&(\lambda_f(n_1)n_1^{-iT}(\log n_1)^m,\lambda_f(n_2)n_2^{-iT}(\log n_2)^m), \\ 
 &(\lambda_f(n_1)n_1^{-i\max\{n_1,n_2\}}(\log n_1)^m, \lambda_f(n_2)n_2^{-i\max\{n_1,n_2\}}(\log n_2)^m)
\end{align*}
in \cite[p.348, LEMMA 6]{GD1}, 
we find that the second and third terms on the right-hand side of \eqref{L2I} are $=O(1)$ uniformly for $\sigma>1$. Hence, the mean square of $L_f^{(m)}(s)$ for ${\rm Re\;}s>1$ is obtained as
\begin{align}
\int_1^T|L_f^{(m)}(\sigma+it)|^2dt=(T-1)\sum_{n=1}^\infty\frac{|\lambda_f(n)|^2(\log n)^{2m}}{n^{2\sigma}}+O(1). \label{LAB}
\end{align}
Next the mean square of $L_f^{(m)}(s)$ for $1/2<{\rm Re\;}s\leq1$ is obtained as follows:
\begin{lem}[{\cite[Theorem 1.3]{YY3}}]
For any $m\in\mathbb{Z}_{\geq0}$ and $T>0$, we have 
\begin{align}
&\int_1^T|L_f^{(m)}(\sigma+it)|^2dt\notag\\
&=\begin{cases} 
\displaystyle (T-1)\sum_{n=1}^\infty\frac{|\lambda_f(n)|^2(\log n)^{2m}}{n^{2\sigma}}+O(T^{2(1-\sigma)}(\log T)^{2m}), & 1/2<\sigma<1, \\ 
\displaystyle (T-1)\sum_{n=1}^\infty\frac{|\lambda_f(n)|^2(\log n)^{2m}}{n^{2\sigma}}+O((\log T)^{2m+2}), & \sigma=1. 
\end{cases}\label{LAC}
\end{align}
\end{lem}
Using Rankin's result mentioned in Introduction
and the following fact
\begin{align*}
\int_{n_f}^\infty\frac{(\log u)^{2m}}{u^{2\sigma}}du=\frac{(2m)!n_f^{1-2\sigma}}{(2\sigma-1)^{2m+1}}\sum_{j=0}^{2m}\frac{(\log{n_f})^j(2\sigma-1)^j}{j!},
\end{align*}
which is obtained by induction, we find that 
the series of main term of (\ref{LAC}) is approximated as
\begin{align}
&\sum_{n=1}^\infty\frac{|\lambda_f(n)|^2(\log n)^{2m}}{n^{2\sigma}}\notag\\
&=-\int_{n_f}^\infty\left(\frac{(\log u)^{2m}}{u^{2\sigma}}\right)'\sum_{n_f<n\leq u}|\lambda_f(n)|^2du\notag\\
&=-\frac{C_f(\log n_f)^{2m}}{n_f^{2\sigma-1}}+C_f\int_{n_f}^{\infty}\frac{(\log u)^{2m}}{u^{2\sigma}}du+O\left(\int_{n_f}^\infty\frac{(\log u)^{2m}}{u^{2\sigma+2/5}}du\right)\notag\\
&=\frac{(2m)!n_fC_f}{n_f^{2\sigma}}\frac{1}{(2\sigma-1)^{2m+1}}+O\left(\frac{1}{(2\sigma-1)^{2m}}\right) \label{LAD}
\end{align}
as $\sigma\to1/2+0$. From (\ref{LAA})--(\ref{LAD}), the following approximate formula is obtained: 
\begin{align}
&\hspace{-12pt}\sum_{\begin{subarray}{c}F(\rho)=0, \\ \sigma\leq{\rm Re\;}\rho\leq\sigma_{f,m},\\ 1\leq{\rm Im\;}\leq T \end{subarray}}({\rm Re\;}\rho-\sigma)\notag\\
\leq&
\displaystyle\frac{(2m+1)(T-1)}{4\pi}\log\frac{1}{2\sigma-1}+\frac{T-1}{4\pi}\log\frac{(2m)!n_fC_f}{|\lambda_f(n_f)|^2(\log n_f)^{2m}}+\notag\\
&+O(\log T)+\frac{T-1}{4\pi}\times \notag\\
&\times\begin{cases} \displaystyle\log\left(1+O\left(\frac{(2\sigma-1)^{2m+1}(\log T)^{2m}}{T^{2\sigma-1}}\right)\right),  & 1/2<\sigma<1, \\
\displaystyle\log\left(1+O\left(\frac{(2\sigma-1)^{2m+1}(\log T)^{2m+2}}{T}\right)\right), \hspace{-6pt}  & \sigma=1,  \\
\displaystyle\log\left(1+O\left(\frac{(2\sigma-1)^{2m+1}}{T}\right)\right), & \sigma{>1}. \end{cases}
\label{LBA}
\end{align}

Finally, we shall give an upper bound of $N_{f,m}(\sigma,T)$. 
Since $N_{f,m}(\sigma,T)$ is monotonically decreasing function with respect to $\sigma$, it follows that
\begin{align}
N_{f,m}(\sigma,T)=&N_{f,m}(\sigma,T)-N_{f,m}(\sigma,1)+C\notag\\
\leq& \frac{1}{\sigma-\sigma_1}\int_{\sigma_1}^{\sigma_{f,m}}(N_{f,m}(u,T)-N_{f,m}(u,1))du+C, \label{LBB}
\end{align}
where we put $\sigma_1=1/2+(\sigma-1/2)/2$. Note that $\sigma-\sigma_1=(\sigma-1/2)/2$, $2\sigma_1-1=\sigma-1/2$.
Since the numbers of zeros of $F_m(s)$ is equal to that of $L_f^{(m)}(s)$, it follows that
\begin{align}
&\int_{\sigma_1}^{\sigma_{f,m}}(N_{f,m}(u,T)-N_{f,m}(u,1))du\notag\\
&=\int_{\sigma_1}^{\sigma_{f,m}}\sum_{\begin{subarray}{c}  F(\rho)=0,\\  u\leq{\rm Re\;}\rho\leq\sigma_{f,m},\\ 1\leq{\rm Im\;}\leq T \end{subarray}}1du
=\sum_{\begin{subarray}{c}  F(\rho)=0,\\  \sigma_1\leq{\rm Re\;}\rho\leq\sigma_{f,m},\\ 1\leq{\rm Im\;}\rho\leq T \end{subarray}}\int_{\sigma_1}^{{\rm Re\;}\rho}1du\notag\\
&=\sum_{\begin{subarray}{c}  F(\rho)=0,\\  \sigma_1\leq{\rm Re\;}\rho\leq\sigma_{f,m},\\ 1\leq{\rm Im\;}\rho\leq T \end{subarray}}({\rm Re\;}\rho-\sigma_1). \label{LBC}
\end{align}
Combining (\ref{LBA})--(\ref{LBC}) we obtain (\ref{ZD1}) and (\ref{ZD2}). 
Hence the proof of Theorem \ref{THM2} is completed.


\bibliographystyle{alpha}

\end{document}